\documentclass{article}
\usepackage{amsmath}
\usepackage{url}
\usepackage{graphicx}
\usepackage{pstricks}
\usepackage{amsfonts}

\begin{document}

\title{Real Analysis in Reverse}
\author{James Propp}
\date{}
\maketitle

\begin{abstract}
Many of the theorems of real analysis,
against the background of the ordered field axioms,
are equivalent to Dedekind completeness,
and hence can serve as completeness axioms for the reals.
In the course of demonstrating this,
the article offers a tour
of some less-familiar ordered fields,
provides some of the relevant history,
and considers pedagogical implications.
\end{abstract}

\newcommand{\bN}{{\mathbb{N}}}
\newcommand{\bZ}{{\mathbb{Z}}}
\newcommand{\bQ}{{\mathbb{Q}}}
\newcommand{\bR}{{\mathbb{R}}}
\newcommand{\bNo}{{\bf No}}

\begin{center}
Copyright the Mathematical Association of America 2012. All rights reserved.
\end{center}

\section{Introduction}

Every mathematician learns that the axioms of an ordered field
(the properties of the numbers 0 and 1, the binary relation $<$, 
and the binary operations $+$, $-$, $\times$, $\div$)
don't suffice as a basis for real analysis;
some sort of heavier-duty axiom is required.
Unlike the axioms of an ordered field,
which involve only % (universal) 
quantification over elements of the field,
the heavy-duty axioms require % (alternating) 
quantification over more complicated objects
such as nonempty bounded {\em subsets} of the field,
Cauchy {\em sequences} of elements of the field,
continuous {\em functions} from the field to itself, 
ways of cutting the field into ``left'' and ``right'' {\em components}, etc.

Many authors of treatises on real analysis remark upon (and prove)
the equivalence of various different axiomatic developments of the theory;
for instance, Korner \cite{refKa} shows that 
the Dedekind Completeness Property
(every nonempty set that is bounded above has a least upper bound)
is equivalent to the Bolzano-Weierstrass Theorem
in the presence of the ordered field axioms.
There are also a number of essays, 
such as Hall's \cite{refHa} and Hall and Todorov's \cite{refHT}, 
that focus on establishing the equivalence of several axioms,
each of which asserts in its own way
that the real number line has no holes in it.
Inasmuch as one of these axioms is the Dedekind Completeness Property,
we call such axioms {\bf completeness properties} for the reals.
(In this article, ``complete'' will always mean ``Dedekind complete'',
except in subsection 5.1.  Readers should be warned that
other authors use ``complete'' to mean ``Cauchy complete''.)
More recently, Teismann \cite{refT} has written an article
very similar to this one, with overlapping aims but differing emphases,
building on an unpublished manuscript by the author \cite{refP}.

One purpose of the current article is to stress that,
to a much greater extent than is commonly recognized,
many theorems of real analysis are completeness properties.
The process of developing this observation is in some ways akin 
to the enterprise of ``reverse mathematics'' 
pioneered by Harvey Friedman and Stephen Simpson
(see e.g.\ \cite{refSd}), wherein one deduces axioms from theorems
instead of the other way around.
However, the methods and aims are rather different.
Reverse mathematics avoids the unrestricted use of ordinary set theory
and replaces it by something tamer, namely, second order arithmetic,
or rather various sub-systems of second order arithmetic
(and part of the richness of reverse mathematics 
arises from the fact that it can matter very much 
which subsystem of second order arithmetic one uses).
In this article, following the tradition of Halmos' classic text 
``Naive Set Theory'' \cite{refHb},
we engage in what might be called naive reverse mathematics,
where we blithely quantify over all kinds of infinite sets
without worrying about what our universe of sets looks like.

Why might a non-logician care about reverse mathematics at all 
(naive or otherwise)?
One reason is that it sheds light on the landscape
of mathematical theories and structures.
Arguably the oldest form of mathematics in reverse
is the centuries-old attempt to determine 
which theorems of Euclidean geometry
are equivalent to the parallel postulate
(see \cite[pp. 276--280]{refGM} for a list of such theorems).
The philosophical import of this work
might be summarized informally as
``Anything that isn't Euclidean geometry 
is very different from Euclidean geometry.''
In a similar way, the main theme of this essay
is that anything that isn't the real number system
must be different from the real number system in many ways.
Speaking metaphoricaly,
one might say that, in the landscape of mathematical theories,
real analysis is an isolated point;
or, switching metaphors, one might say that 
the real number system is rigid
in the sense that it cannot be subjected to slight deformations.

An entertaining feature of real analysis in reverse
is that it doesn't merely show us how 
some theorems of the calculus that look different
are in a sense equivalent;
it also shows us how some theorems that look fairly similar
are {\em not} equivalent.
Consider for instance the following three propositions
taught in calculus classes:

\noindent
(A) The Alternating Series Test:
If $a_1 \geq a_2 \geq a_3 \geq \dots$
and $a_n \rightarrow 0$,
then $\sum_{n=1}^{\infty} (-1)^n a_n$ converges.

\noindent
(B) The Absolute Convergence Theorem: 
If $\sum_{n=1}^\infty |a_n|$ converges,
then $\sum_{n=1}^\infty a_n$ converges.

\noindent
(C) The Ratio Test: If $|a_{n+1}/a_n| \rightarrow L$ as $n \rightarrow \infty$,
with $L < 1$, then $\sum_{n=1}^{\infty} a_n$ converges.

\noindent
To our students (and perhaps to ourselves)
the three results can seem much of a muchness,
yet there is a sense in which one of the three theorems
is stronger than the other two.
Specifically, one and only one of them is equivalent to completeness
(and therefore implies the other two).
How quickly can you figure out which is the odd one out?

At this point some readers may be wondering what I mean by equivalence.
(``If two propositions are theorems, don't they automatically imply each other,
since by the rules of logic
every true proposition materially implies every other true proposition?'')
Every proposition $P$ in real analysis,
being an assertion about $\bR$,
can be viewed more broadly as a family of assertions $P(R)$
about ordered fields $R$;
one simply takes each explicit or implicit reference to $\bR$
in the proposition $P$
and replaces it by a reference to some unspecified ordered field $R$.
Thus, every theorem $P$ 
is associated with a property $P(\cdot)$ satisfied by $\bR$
and possibly other ordered fields as well.
What we mean when we say that one proposition of real analysis $P$
{\it implies\/} another proposition of real analysis $P'$
is that $P'(R)$ holds whenever $P(R)$ holds,
where $R$ varies over all ordered fields;
and what we mean by the {\it equivalence\/} of $P$ and $P'$
is that this relation holds in both directions.
In particular, when we say that $P$ is a completeness property,
or that it can serve as an axiom of completeness,
what we mean is that for any ordered field $R$,
$P(R)$ holds if and only if $R$ satisfies Dedekind completeness.
(In fact, Dedekind proved \cite[p.\ 33]{refBr}
that any two ordered fields that are Dedekind complete are isomorphic;
that is, the axioms of a Dedekind complete ordered field are {\bf categorical}.)

To prove that a property $P$ satisfied by the real numbers
is {\it not\/} equivalent to completeness,
we need to show that there exists an ordered field
that satisfies property $P$ but not the completeness property. 
So it's very useful to have on hand a number of different ordered fields
that are {\it almost\/} the real numbers, but not quite.
The second major purpose of the current article
is to introduce the reader to some ordered fields of this kind.
We will often call their elements ``numbers'',
since they behave like numbers in many ways.
(This extension of the word ``number'' is standard
when one speaks of a different variant of real numbers, 
namely $p$-adic numbers;
however, this article is about {\em ordered} fields,
so we will have nothing to say about $p$-adics.)

A third purpose of this article is 
to bring attention to Dedekind's Cut Property
(property (3) of section 2).
Dedekind singled out this property of the real numbers
as encapsulating what makes $\bR$ a continuum,
and if the history of mathematics had gone slightly differently,
this principle would be part of the standard approach to the subject.
However, Dedekind never used this property as the basis
of an axiomatic approach to the real numbers;
instead, he constructed the real numbers from the rational numbers
via Dedekind cuts and then verified that the Cut Property holds.
Subsequently, most writers of treatises and textbooks on real analysis
adopted the Least Upper Bound Property 
(aka the Dedekind Completeness Property) as the
heavy-duty second order axiom that distinguishes 
the real number system from its near kin.
And indeed the Least Upper Bound Property 
is more efficient than the Cut Property
for purposes of getting the theory of the calculus off the ground.
But the Cut Property has a high measure of symmetry and simplicity
that is missing from its rival.
You can explain it to average calculus students,
and even lead them to conjecture it on their own;
the only thing that's hard is convincing them
that it's nontrivial!
The Cut Property hasn't been entirely forgotten
(\cite{refAS}, \cite[p.\ 53]{refKb}, and \cite{refMe})
and it's well-known among people who study
the axiomatization of Euclidean geometry \cite{refG}
or the theory of partially ordered sets and lattices \cite{refWaa}.
But it deserves to be better known
among the mathematical community at large.

This brings me to my fourth purpose, which is pedagogical.
There is an argument to be made that, in the name of intellectual honesty, 
we who teach more rigorous calculus courses (often billed as ``honors'' courses)
should try to make it clear what assumptions
the theorems of the calculus depend on,
even if we skip some (or most) of the proofs
in the chain of reasoning that leads
from the assumptions to the central theorems of the subject,
and even if the importance of the assumptions will not be fully clear
until the students have taken more advanced courses.
It is most common to use the Dedekind Completeness Property 
or Monotone Sequence Convergence Property for this purpose,
and to introduce it explicitly only late in the course,
after differentiation and integration have been studied,
when the subject shifts to infinite sequences and series.
I will suggest some underused alternatives.

Note that this article is {\it not} about ways of constructing the real numbers.
(The Wikipedia page \cite{refWa} gives 
both well-known constructions and more obscure ones;
undoubtedly many others have been proposed.)
This article is about the axiomatic approach to real analysis,
and the ways in which the real number system
can be characterized by its internal properties.

The non-introductory sections of this article are structured as follows.
In the second section,
I'll state some properties of ordered fields $R$ that hold when $R = \bR$.
In the third section, I'll give some examples of ordered fields
that resemble (but aren't isomorphic to) the field of real numbers.
In the fourth section,
I'll show which of the properties from the second section
are equivalent to (Dedekind) completeness and which aren't.
In the fifth section, I'll discuss some tentative pedagogical implications.

It's been fun for me to write this article,
and I have a concrete suggestion for how
the reader may have the same kind of fun.
Start by reading the second section
and trying to decide on your own which of the properties
are equivalent to completeness and which aren't.
If you're stumped in some cases,
and suspect but cannot prove
that some particular properties aren't equivalent to completeness,
read the third section
to see if any of the ordered fields discussed there
provide a way to see the inequivalence.
And if you're still stumped,
read the fourth section.
You can treat the article as a collection of math puzzles,
some (but not all) of which can be profitably contemplated in your head.

A reminder about the ground-rules for these puzzles is in order.
Remember that we are to interpret every theorem of real analysis
as the particular case $R=\bR$
of a family of propositions $P(R)$ about ordered fields $R$.
An ordered field is a collection of elements
(two of which are named $0_R$ and $1_R$,
or just $0$ and $1$ for short),
equipped with the relations $<$
and operations $+$, $-$, $\times$, $\div$
satisfying all the usual ``high school math'' properties.
(Note that by including subtraction and division as primitives,
we have removed the need for existential quantifiers in our axioms;
e.g., instead of asserting that for all $x \neq 0$
there exists a $y$ such that $x \times y = y \times x=1$,
we simply assert that for all $x \neq 0$, 
$x \times (1 \div x) = (1 \div x) \times x = 1$.
It can be argued that instead of minimizing 
the number of primitive notions or the number of axioms,
axiomatic presentations of theories should minimize 
the number of existential quantifiers,
and indeed this is standard practice in universal algebra.)
In any ordered field $R$ we can define notions like
$|x|_R$ (the unique $y \geq 0$ in $R$ with $y = \pm x$)
and $(a,b)_R$ (the set of $x$ in $R$ with $a < x < b$).
More complicated notions from real analysis
can be defined as well:
for instance, given a function $f$ from $R$ to $R$,
we can define $f'(a)_R$, the ``derivative of $f$ at $a$ relative to $R$'',
as the $L \in R$ (unique if it exists)
such that for every $\epsilon > 0_R$
there exists $\delta > 0_R$
such that $|(f(x)-f(a))/(x-a) - L|_R < \epsilon$
for all $x$ in $(a-\delta,a)_R \cup (a,a+\delta)_R$.
The subscripts are distracting, so we will omit them,
but it should be borne in mind that they are conceptually present.
What goes for derivatives of functions
applies to other notions of real analysis as well,
such as the notion of a convergent sequence:
the qualifier ``relative to $R$'', even if unstated,
should always be kept in mind.

The version of naive set theory we will use
includes the axiom of countable choice.
A reasonably large subset of real analysis can be set up
without countable choice,
but many important theorems such as Bolzano-Weierstrass
rely on countable choice in an essential way.
Whether or not one believes that the axiom of countable choice is true,
distrusting countable choice
requires a fair amount of foundational sophistication,
and therefore cannot in my opinion
be considered a truly ``naive'' stance.

Every ordered field $R$ contains an abelian semigroup 
$\bN_R = \{ 1_R, \ 1_R+1_R, \ 1_R+1_R+1_R, \ \dots\}$ isomorphic to $\bN$;
$\bN_R$ may be described as the intersection of all subsets of $R$
that contain $1_R$ and are closed under the operation
that sends $x$ to $x+1_R$.
Likewise every ordered field $R$ contains
an abelian group $\bZ_R$ isomorphic to $\bZ$,
and a subfield $\bQ_R$ isomorphic to $\bQ$.
We shall endow an ordered field $R$ with the {\bf order topology},
that is, the topology generated by basic open sets of the form $(a,b)_R$.

\section{Some theorems of analysis}

The following propositions about an ordered field $R$
(and about associated structures such as 
$[a,b]=[a,b]_R=\{x \in R:\ a \leq x \leq b\}$)
are true when the ordered field $R$
is taken to be $\bR$, the field of real numbers.

(1) The Dedekind Completeness Property:
Suppose $S$ is a nonempty subset of $R$ that is bounded above.
Then there exists a number $c$ that is an upper bound of $S$
such that every upper bound of $S$ is greater than or equal to $c$.

(2) The Archimedean Property: 
For every $x \in R$ there exists $n \in \bN_R$ with $n > x$.
Equivalently, for every $x \in R$ with $x > 0$
there exists $n \in \bN_R$ with $1/n < x$.

(3) The Cut Property: Suppose $A$ and $B$ are nonempty disjoint subsets of $R$
whose union is all of $R$, such that every element of $A$ is less than
every element of $B$.  Then there exists a cutpoint $c \in R$ such that
every $x < c$ is in $A$ and every $x>c$ is in $B$.
(Or, if you prefer: Every $x \in A$ is $\leq c$,
and every $x \in B$ is $\geq c$.
It's easy to check that the two versions are equivalent.)
Since this property may be unfamiliar, we remark that
the Cut Property follows immediately from Dedekind completeness
(take $c$ to be the least upper bound of $A$).
 
(4) Topological Connectedness: 
Say $S \subseteq R$ is open if for every $x$ in $S$
there exists $\epsilon > 0$ so that
every $y$ with $|y-x| < \epsilon$ is also in $S$.
Then there is no way to express $R$
as a union of two disjoint nonempty open sets.
That is, if $R = A \cup B$ with $A,B$ nonempty and open,
then $A \cap B$ is nonempty.

(5) The Intermediate Value Property: 
If $f$ is a continuous function
from $[a,b]$ to $R$, with $f(a) < 0$ and $f(b) > 0$,
then there exists $c$ in $(a,b)$ with $f(c) = 0$.

(6) The Bounded Value Property:
If $f$ is a continuous function from $[a,b]$ to $R$,
there exists $B$ in $R$ with $f(x) \leq B$ for all $x$ in $[a,b]$.

(7) The Extreme Value Property:
If $f$ is a continuous function from $[a,b]$ to $R$,
there exists $c$ in $[a,b]$ with $f(x) \leq f(c)$ for all $x$ in $[a,b]$.

(8) The Mean Value Property: 
Suppose $f: [a,b] \rightarrow R$ is continuous on $[a,b]$ 
and differentiable on $(a,b)$.
Then there exists $c$ in $(a,b)$ such that
$f'(c) = (f(b)-f(a))/(b-a)$.

(9) The Constant Value Property:
Suppose $f: [a,b] \rightarrow R$ is continuous on $[a,b]$ 
and differentiable on $(a,b)$,
with $f'(x) = 0$ for all $x$ in $(a,b)$.
Then $f(x)$ is constant on $[a,b]$.

(10) The Convergence of Bounded Monotone Sequences:
Every monotone increasing (or decreasing) sequence in $R$
that is bounded converges to some limit.

(11) The Convergence of Cauchy Sequences: 
Every Cauchy sequence in $R$ is convergent.

(12) The Fixed Point Property for Closed Bounded Intervals:
Suppose $f$ is a continuous map from $[a,b] \subset R$ to itself.
Then there exists $x$ in $[a,b]$ such that $f(x) = x$.

(13) The Contraction Map Property: 
Suppose $f$ is a map from $R$ to itself such that for some constant $c < 1$,
$|f(x) - f(y)| \leq c|x-y|$ for all $x,y$.
Then there exists $x$ in $R$ such that $f(x) = x$.

(14) The Alternating Series Test: 
If $a_1 \geq a_2 \geq a_3 \geq \dots$
and $a_n \rightarrow 0$,
then $\sum_{n=1}^{\infty} (-1)^n a_n$ converges.

(15) The Absolute Convergence Property: 
If $\sum_{n=1}^\infty |a_n|$ converges in $R$,
then $\sum_{n=1}^\infty a_n$ converges in $R$.

(16) The Ratio Test: 
If $|a_{n+1}/a_n| \rightarrow L$ in $R$ as $n \rightarrow \infty$,
with $L < 1$, then $\sum_{n=1}^{\infty} a_n$ converges in $R$.

(17) The Shrinking Interval Property:
Suppose $I_1 \supseteq I_2 \supseteq \dots$ 
are bounded closed intervals in $R$ with lengths decreasing to 0.
Then the intersection of the $I_n$'s is nonempty.

(18) The Nested Interval Property:
Suppose $I_1 \supseteq I_2 \supseteq \dots$ 
are bounded closed intervals in $R$.
Then the intersection of the $I_n$'s is nonempty.
% (Another name for this property is ``spherical completeness''.)

\section{Some ordered fields}

The categoricity of the axioms for $\bR$
tells us that any ordered field that is Dedekind-complete
must be isomorphic to $\bR$.
So one (slightly roundabout) way to see that
the ordered field of rational numbers $\bQ$ fails to satisfy completeness
is to note that it contains too few numbers to be isomorphic to $\bR$.
The same goes for the field of real algebraic numbers.
There are even bigger proper subfields of the $\bR$;
for instance, Zorn's Lemma implies that 
among the ordered subfields of $\bR$ that don't contain $\pi$,
there exists one that is maximal with respect to this property.
% Consider a totally ordered collection of subfields of R.
% Call the subfields F_i, where i ranges over some index set I.
% Suppose each F_i is a subfield of R not containing R. 
% Let F be the union of the F_i's.
% It's a subfield of R not containing pi.
But most of the ordered fields
we'll wish to consider have the opposite problem:
they contain too {\em many} ``numbers''.

Such fields may be unfamiliar, but logic tells us
that number systems of this kind must exist.
(Readers averse to ``theological'' arguments might prefer
to skip this paragraph and the next
and proceed directly to a concrete construction of such a number system
two paragraphs below,
but I think there is value in an approach that convinces us ahead of time 
that the goal we seek is not an illusory one,
and shows that such number systems exist
without commiting to one such system in particular.)
Take the real numbers and adjoin 
a new number $n$, satisfying the infinitely many axioms $n > 1$,
$n > 2$, $n > 3$, etc.  
Every finite subset of this infinite set of first-order axioms 
(together with the set of ordered-field axioms)
has a model, so by the compactness principle of first-order logic
(see e.g.\ \cite{refDM}),
these infinitely many axioms must have a model.  
(Indeed, I propose that the compactness principle
is the core of validity inside
the widespread student misconception that 0.999\dots is different from 1; 
on some level, students may be reasoning that if the intervals 
$[0.9,1.0)$, $[0.99,1.00)$, $[0.999,1.000)$, etc.\ are all nonempty, 
then their intersection is nonempty as well.
The compactness principle tells us that
there must exist ordered fields in which
the intersection of these intervals is nonempty. 
Perhaps we should give these students credit for intuiting, in a murky way,
the existence of non-Archimedean ordered fields!)

But what does an ordered field with infinite elements
(and their infinitesimal reciprocals) look like?

One such model is given by rational functions in one variable,
ordered by their behavior at infinity;
we call this variable $\omega$ rather than the customary $x$,
since it will turn out to be bigger than every real number,
under the natural imbedding of $\bR$ in $R$.
Given two rational functions $q(\omega)$ and $q'(\omega)$
decree that $q(\omega) > q'(\omega)$ iff $q(r) > q'(r)$
for all sufficiently large real numbers $r$.
One can show (see e.g.\ \cite{refKa})
that this turns $\bR(\omega)$ into an ordered field.
We may think of our construction as the process of
adjoining a formal infinity to $\bR$.
Alternatively, we can construct an ordered field isomorphic to $\bR(\omega)$
by adjoining a formal infinitesimal $\epsilon$
(which the isomorphism identifies with $1/\omega$):
given two rational functions $q(\varepsilon)$ and $q'(\varepsilon)$,
decree that $q(\varepsilon) > q'(\varepsilon)$ iff $q(r) > q'(r)$
for all positive real numbers $r$ sufficiently close to 0.
Note that this ordered field is non-Archimedean:
just as $\omega$ is bigger than every real number,
the positive element $\varepsilon$
is less than every positive real number.

We turn next to the field of formal Laurent series.
A formal Laurent series is a formal expression
$\sum_{n \geq N} a_n \varepsilon^n$
where $N$ is some non-positive integer
and the $a_n$'s are arbitrary real numbers;
the associated finite sum $\sum_{N \leq n < 0} a_n \varepsilon^n$
is called its {\bf principal part}.
One can define field operations on such expressions
by mimicking the ordinary rules of adding, subtracting,
multiplying, and dividing Laurent series,
without regard to issues of convergence.
The leading term of such an expression is the nonvanishing term
$a_n \varepsilon^n$ for which $n$ is as small as possible,
and we call the expression positive or negative
according to the sign of the leading term.
In this way we obtain an ordered field.
This field is denoted by $\bR((\varepsilon))$,
and the field $\bR(\varepsilon)$ discussed above
may be identified with a subfield of it.
In this larger field, a sequence of Laurent series converges
if and only if the sequence of principal parts stabilizes
(i.e., is eventually constant)
and for every integer $n \geq 0$
the sequence of coefficients of $\varepsilon^n$ stabilizes.
In particular, $1,\varepsilon,\varepsilon^2,\varepsilon^3,\dots$
converges to 0
but $1,\frac12,\frac14,\frac18,\dots$ does not.
The same holds for $R = \bR(\varepsilon)$;
the sequence $1,\frac12,\frac14,\frac18,\dots$ 
does not converge to 0 relative to $R$
because every term differs from 0 by more than $\varepsilon$.

Then come the really large non-Archimedean ordered fields.
There are non-Archimedean ordered fields so large
(that is, equipped with so many infinitesimal elements)
that the ordinary notion of convergence of sequences becomes trivial:
all convergent sequences are eventually constant.
In such an ordered field,
the only way to ``sneak up'' on an element from above or below
is with a generalized sequence
whose terms are indexed by some uncountable ordinal,
rather than the countable ordinal $\omega$.
Define the {\bf cofinality} of an ordered field
as the smallest possible cardinality of
an unbounded subset of the field
(so that for instance the real numbers,
although uncountable, have countable cofinality).
The cardinality and cofinality of a non-Archimedean ordered field
can be as large as you like (or dislike!);
along with Cantor's hierarchies of infinities
comes an even more complicated hierarchy of non-Archimedean ordered fields.
The cofinality of an ordered field is easily shown to be a regular cardinal, 
%% check this!
where a cardinal $\kappa$ is called {\bf regular} iff 
a set of cardinality $\kappa$ cannot be written as the union 
of fewer than $\kappa$ sets each of cardinality less than $\kappa$.
In an ordered field $R$ of cofinality $\kappa$,
the right notion of a sequence is a ``$\kappa$-sequence'',
defined as a function from the ordinal $\kappa$
(that is, from the set of all ordinals $\alpha < \kappa$) to $R$.
% (Here I am using $\kappa$ to denote both an ordinal and a cardinal;
% this double usage is standard and usually does not lead to confusion.)
A sequence whose length is less than the cofinality of $R$
can converge only if it is eventually constant.
% Equivalently, a sequence whose length is less than the cofinality of R
% can converge to 0 only if it is eventually 0.
% Proof: Suppose (x_\alpha : \alpha < \kappa') 
% is a sequence of length \kappa' < \kappa.
% To say the sequence converges to 0 is to say that
% for any epsilon, the sequence eventually lies in (-epsilon,+epsilon).
% If it is not eventually constant, this means that there must be
% at least one nonzero term of the sequence in (-epsilon,+epsilon).
% That is, there exists a nonzero term x_\alpha with |x_\alpha| < epsilon,
% so that |1/x_\alpha| > 1/epsilon.  I.e., for every B ( = 1/epsilon),
% there is an element x of the set S = {|1/x_\alpha|: x_\alpha \neq 0}
% with x > B.  Hence S is unbounded.  But S has cardinality kappa'.

Curiously, if one uses this generalized notion of a sequence,
some of the large ordered fields
can be seen to have properties of generalized compactness
reminiscent of the real numbers.
More specifically, for $\kappa$ regular,
say that an ordered field $R$ of cofinality $\kappa$
satisfies the $\kappa$-Bolzano-Weierstrass Property
if every bounded $\kappa$-sequence $(x_{\alpha})_{\alpha<\kappa}$ in $R$
has a convergent $\kappa$-subsequence.
Then a theorem of Sikorski \cite{refSc}
states that for every uncountable regular cardinal $\kappa$
there is an ordered field of cardinality and cofinality $\kappa$ 
that satisfies the $\kappa$-Bolzano-Weierstrass Property.
For more background on non-Archimedean ordered fields
(and generalizations of the Bolzano-Weierstrass Property in particular),
readers can consult \cite{refJS} and/or \cite{refSa}.

Lastly, there is the Field of surreal numbers $\bNo$,
which contains {\it all\/} ordered fields as subfields.
Following Conway \cite{refCo} we call it a Field rather than a field
because its elements form a proper class rather than a set.
One distinguishing property of the surreal numbers
is the fact that for any two sets of surreal numbers $A,B$
such that every element of $A$ is less than every element of $B$,
there exists a surreal number
that is greater than every element of $A$
and less than every element of $B$.
(This does not apply if $A$ and $B$ are proper classes,
as we can see by letting $A$ consist of 0 and the negative surreal numbers
and $B$ consist of the positive surreal numbers.)

See the Wikipedia page \cite{refWb} for information on other ordered fields,
such as the Levi-Civita field and the field of hyper-real numbers.

\section{Some proofs}

Here we give (sometimes abbreviated) versions of the proofs
of equivalence and inequivalence.

\bigskip

The Archimedean Property (2) 
does not imply the Dedekind Completeness Property (1):
The ordered field of rational numbers satisfies the former but not the latter.
(Note however that (1) does imply (2):
$\bN_R$ is nonempty, so if $\bN_R$ were bounded above, 
it would have a least upper bound $c$ by (1).
Then for every $n \in \bN_R$ we would have $n+1 \leq c$
(since $n+1$ is in $\bN_R$ and $c$ is an upper bound for $\bN_R$),
implying $n \leq c-1$.
But then $c-1$ would be an upper bound for $\bN_R$,
contradicting our choice of $c$ as a least upper bound for $\bN_R$.
This shows that $\bN_R$ is not bounded above, which is (2).)

\bigskip

The Cut Property (3) implies completeness (1):
Given a nonempty set $S \subseteq R$ that is bounded above,
let $B$ be the set of upper bounds of $S$
and $A$ be its complement.
$A$ and $B$ satisfy the hypotheses of (3),
so there exists a number $c$ such that
everything less than $c$ is in $A$ and everything greater than $c$ is in $B$.
It is easy to check that $c$ is a least upper bound of $S$.
(To show that $c$ is an upper bound of $S$,
suppose some $s$ in $S$ exceeds $c$.
Since $(s+c)/2$ exceeds $c$, it belongs to $B$,
so by the definition of $B$ it must be an upper bound of $S$,
which is impossible since $s > (s+c)/2$.
To show that $c$ is a least upper bound of $S$,
suppose that some $a < c$ is an upper bound of $S$.
But $a$ (being less than $c$) is in $A$, 
so it can't be an upper bound of $S$.)

\bigskip

In view of the preceding result, the Cut Property
is a completeness property,
and to prove that some other property is a completeness property,
it suffices to show that it implies the Cut Property.
Hereafter we will write
``Property ($n$) implies completeness by way of the Cut Property (3)''
to mean ``($n$) $\Rightarrow$ (3) $\Rightarrow$ (1).''
When a detour through the Archimedean Property is required
as a lemma to the proof of the Cut Property, we will write
``Property ($n$) implies completeness by way of 
the Archimedean Property (2) and the Cut Property (3)''
to give a road-map of the argument that follows.

\bigskip

Topological Connectedness (4) implies completeness 
by way of the Cut Property (3):
We prove the contrapositive.
Let $A$ and $B$ be sets satisfying the hypotheses of the Cut Property
but violating its conclusion:
there exists no $c$ such that
everything less than $c$ is in $A$ and everything greater than $c$ is in $B$.
(That is, suppose $A,B$ is a ``bad cut'', which we also call a {\bf gap}.)
Then for every $a$ in $A$ there exists $a'$ in $A$ with $a' > a$,
and for every $b$ in $B$ there exists $b'$ in $B$ with $b' < b$.
From this it readily follows that the sets $A$ and $B$ are open, 
so that Topological Connectedness must fail.

\bigskip

The Intermediate Value Property (5) implies completeness 
by way of the Cut Property (3):
We again prove the contrapositive.
(Indeed, we will use this mode of proof so often
that henceforth we will omit the preceding prefatory sentence.)
Let $A,B$ be a gap.
The function that is $-1$ on $A$ and $1$ on $B$
is continuous and violates the conclusion of the Intermediate Value Property.

\bigskip

The Bounded Value Property (6) does not imply completeness:
Counterexamples are provided by the theorem of Sikorski
referred to earlier (near the end of Section 3),
once we prove that the $\kappa$-Bolzano-Weirestrass Property 
implies the Bounded Value Property.
(I am indebted to Ali Enayat for suggesting this approach
and for supplying all the details that appear below.)

Suppose that $\kappa$ is a regular cardinal
and that $R$ is an ordered field with cofinality $\kappa$
satisfying the $\kappa$-Bolzano-Weierstrass Property.
Then I claim that $R$ satisfies the Bounded Value Property.
For, choose an increasing unbounded sequence 
$(x_\alpha \, : \, \alpha \in \kappa)$ of elements of $R$.
Suppose that $f$ is continuous on $[a,b]$
but that there exists no $B$ with $f(x) \leq B$ for all $x \in [a,b]$.
For each $\alpha \in \kappa$,
choose some $t_\alpha \in [a,b]$ with $f(t_{\alpha}) > x_{\alpha}$.
The $t_\alpha$'s are bounded, so by the $\kappa$-Bolzano-Weierstrass Property
there exists some subset $U$ of $\kappa$
such that $(t_\alpha \, : \, \alpha \in U)$ is a $\kappa$-subsequence 
that converges to some $c \in [a,b]$.
(Note that $U$ must be unbounded;
otherwise $(t_\alpha \, : \, \alpha \in U)$
would be a $\beta$-sequence for some $\beta < \kappa$
rather than a $\kappa$-sequence.)
By the continuity of $f$, 
the sequence $(f(t_\alpha) \, : \, \alpha \in U)$ converges to $f(c)$.
% We need to verify that continuous maps respect convergence under
% this transfinite ordinal notion of convergence.  Let's check it:
% We want to show that if x_n converges to x, then f(x_n) converges to f(x).  
% Fix epsilon > 0.  By continuity of f at x, there exists delta > 0 such 
% that every x' with |x'-x| < delta satisfies |f(x')-f(x)| < epsilon.
% Since x_n converges to x, there exists an N such that |x_n - x| < delta
% for all n > N.  But then |f(x_n)-f(x)| < epsilon for all n > N.
% Since epsilon was arbitrary, f(x_n) converges to f(x), as claimed.
% In short, it's the standard proof, but with ordinals.

We now digress to prove a small lemma, namely,
that every convergent $\kappa$-sequence 
$(r_\alpha \in R \, : \, \alpha \in \kappa)$ is bounded.
Since $(r_\alpha \, : \, \alpha \in \kappa)$ 
converges to some limit $r$,
there exists a $\beta < \kappa$ such that for all $\alpha \geq \beta$,
$r_\alpha$ lies in $(r-1,r+1)$;
then the tail-set $\{r_\alpha \, : \, \alpha \geq \beta\}$ is bounded.
On the other hand,
since $\kappa$ is a regular cardinal, and since $R$ has cofinality $\kappa$,
the complementary set
$\{r_\alpha \, : \, \alpha < \beta\}$ is too small to be unbounded.
Taking the union of these two sets,
we see that the set $\{r_\alpha \in R \, : \, \alpha \in \kappa\}$ is bounded.
%% Look this over again later to make sure I understand it.

Applying this lemma to the convergent sequence
$(f(t_\alpha) \, : \, \alpha \in U)$,
we see that $(f(t_\alpha) \, : \, \alpha \in U)$ is bounded.
But this is impossible, since the set $U$ is unbounded
and since our original increasing sequence
$(x_\alpha \, : \, \alpha \in \kappa)$ was unbounded.
This contradiction shows that $f([a,b])$ is bounded above.
Hence $R$ has the Bounded Value Property, as claimed.

(Note the resemblance between the preceding proof
and the usual real-analysis proof
that every continuous real-valued function on an interval $[a,b]$ is bounded.)

One detail omitted from the above argument
is a proof that uncountable regular cardinals exist
(without which Sikorski's theorem is vacuous).
The axiom of countable choice implies
that a countable union of countable sets is countable,
so $\aleph_1$, the first uncountable cardinal, is a regular cardinal.
% However, it would be good to have a way to prove
% that the Bounded Value Property (6) does not imply completeness
% without requiring the Axiom of Choice
% (or else an explanation of why no such proof could exist).
%% Am I being inconsistent here?  Do I use choice elsewhere
%% without commenting upon it, as in the case of other
%% counterexamples involving big ordered fields?

\bigskip

The Extreme Value Property (7) implies completeness 
by way of the Cut Property (3):
Suppose $A,B$ is a gap, and for convenience assume $1 \in A$ and $2 \in B$
(the general case may be obtained from this special case by
straightforward algebraic modifications).  Define 
$$f(x) = \left\{ \begin{array}{ll}
	x & \mbox{if $x \in A$}, \\
	0 & \mbox{if $x \in B$}
\end{array} \right.$$
for $x$ in $[0,3]$.
Then $f$ is continuous on $[0,3]$
but there does not exist $c \in [0,3]$
with $f(x) \leq f(c)$ for all $x$ in $[0,3]$.
For, such a $c$ would have to be in $A$
(since $f$ takes positive values on $[0,3] \cap A$, e.g.\ at $x=1$,
while $f$ vanishes on $[0,3] \cap B$),
but for any $c \in [0,3] \cap A$
there exists $c' \in [0,3] \cap A$ with $c' > c$,
so that $f(c') > f(c)$.

\bigskip

The Mean Value Property (8) implies the Constant Value Property (9): Trivial.

\bigskip

The Constant Value Property (9) implies completeness 
by way of the Cut Property (3):
Suppose $A,B$ is a gap.
Again consider the function $f$ that equals $-1$ on $A$
and 1 on $B$.  It has derivative 0 everywhere, 
yet it isn't constant on $[a,b]$
if one takes $a \in A$ and $b \in B$.

\bigskip

The Convergence of Bounded Monotone Sequences (10)
implies completeness by way of the Archimedean Property (2) 
and the Cut Property (3):
If $R$ does not satisfy the Archimedean Property,
then there must exist an element of $R$
that is greater than every term of the sequence $1,2,3,\dots$.
By the Convergence of Bounded Monotone Sequences,
this sequence must converge, say to $r$.
This implies that $0,1,2,\dots$ also converges to $r$.
Now subtract the two sequences;
by the algebraic limit laws that are easily derived from
the ordered field axioms and the definition of limits,
one finds that $1,1,1\dots$ converges to 0,
which is impossible.
Therefore $R$ must satisfy the Archimedean Property.
Now suppose we are given a cut $A,B$.
For $n \geq 0$ in $\bN$,
let $a_n$ be the largest element of $2^{-n} \bZ_R$ in $A$
and $b_n$ be the smallest element of $2^{-n} \bZ_R$ in $B$.
$(a_n)$ and $(b_n)$ are bounded monotone sequence,
so by the Convergence of Bounded Monotone Sequences they converge, 
and since (by the Archimedean Property) $a_n - b_n$ converges to 0,
$a_n$ and $b_n$ must converge to the same limit; call it $c$.
We have $a_n \leq c \leq b_n$ for all $n$,
so $|a_n - c|$ and $|b_n - c|$ are both at most $2^{-n}$.
From the Archimedean Property it follows that
for every $\epsilon > 0$ there exists $n$ with $2^{-n} < \epsilon$,
and for this $n$ we have $a_n \in A$ and $b_n \in B$
satisfying $a_n > c - \epsilon$ and $b_n < c + \epsilon$
Hence for every $\epsilon > 0$,
every number less than or equal to $c - \epsilon$ is in $A$
and every number greater than or equal to $c + \epsilon$ is in $B$.
Therefore every number less than $c$ is in $A$
and every number greater than $c$ is in $B$,
which verifies the Cut Property.

\bigskip

The Convergence of Cauchy Sequences (11) does not imply completeness:
Every Cauchy sequence in the field of formal Laurent series converges,
but the field does not satisfy the Cut Property.
Call a formal Laurent series {\bf finite} if it is of the form
$\sum_{n \geq 0} a_n \varepsilon^n$;
otherwise, call it {\bf positively infinite} or {\bf negatively infinite}
according to the sign of its leading term $a_n$ ($n < 0$).
If we let $A$ be the set of 
all finite or negatively infinite formal Laurent series,
and we let $B$ be the set of all positively infinite formal Laurent series,
then $A,B$ is a gap.
On the other hand, it is easy to show that
this ordered field satisfies property (11).
Note also that if one defines the norm
of a formal Laurent series $\sum_{n \geq N} a_n \varepsilon^n$
(with $a_N \neq 0$) as $2^{-N}$
and defines the distance between two series
as the norm of their difference,
then one obtains a complete metric space
whose metric topology coincides with
the order topology introduced above.

\bigskip

The Fixed Point Property for Closed Bounded Intervals (12)
implies completeness by way of the Cut Property (3):
Let $A,B$ be a gap of $R$.
Pick $a$ in $A$ and $b$ in $B$, and define $f:[a,b] \rightarrow R$
by putting $f(x) = b$ for $x \in A$ and $f(x) = a$ for $x \in B$.
Then $f$ is continuous but has no fixed point.

\bigskip

The Contraction Map Property (13) implies completeness by way of 
the Archimedean Property (2) and the Cut Property (3):
Here is an adaptation of a solution found by George Lowther \cite{refL}.
First we will show that $R$ is Archimedean.  Suppose not.
Call $x$ in $R$ {\bf finite} if $-n < x < n$ for some $n$ in $\bN_R$,
and {\bf infinite} otherwise.  Let
$$f(x) = \left\{ \begin{array}{ll}
	\frac12 x        & \mbox{if $x$ is infinite}, \\
	x + \frac12 g(x) & \mbox{if $x$ is finite}
\end{array} \right.$$
with 
$$g(x) = 1 - \frac{x}{1+|x|},$$
a decreasing function of $x$ taking values in $(0,2)$.
For all finite $x,y$ with $x>y$ 
one has $(g(y)-g(x))/(x-y) \geq 1/(1+|x|)(1+|y|)$
(indeed, the left-hand side minus the right-hand side 
equals 0 in the cases $x > y \geq 0$ or $0 \geq x > y$
and equals $-2xy/(1+x)(1-y)(x-y) > 0$ in the case $x > 0 > y$),
so for all $x>y$ in $[-a,a]$ we have
$(g(y)-g(x))/(x-y) \geq 1/(1+a)^2$,
implying that $|(f(x)-f(y))/(x-y)| \leq 1 - \frac12 / (1+a)^2$
for all $x,y$ in $[-a,a]$.
Taking $c = 1 - \frac12 / \omega^2$
with $\omega > n$ for all $n$ in $\bN_R$,
we obtain $|f(x)-f(y)| < c|x-y|$
for all finite $x,y$.
This inequality can also be shown to hold 
when one or both of $x,y$ is infinite.
Hence $f$ is a contraction map, 
yet it has no finite or infinite fixed points; contradiction.

Now we want to prove that $R$ satisfies the Cut Property.  Suppose not.
Let $A,B$ be a gap of $R$, and 
let $a_n = \max \, A \cap 2^{-n} \bZ_R$
and $b_n = \min \, B \cap 2^{-n} \bZ_R$.
Since $A,B$ is a gap,
neither of the sequences $(a_n)_{n=1}^{\infty}$ and $(b_n)_{n=1}^{\infty}$
can be eventually constant,
so there exist $n_1 < n_2 < n_3 < \dots$
such that the sequences $(x_k)_{k=1}^\infty$ and $(y_k)_{k=1}^\infty$
with $x_k = a_{n_k}$ and $y_k = b_{n_k}$
are strictly monotone,
with $0 < y_k - x_k \leq 2^{-k}$.
By the Archimedean Property,
every element of $A$ lies in $(-\infty,x_1]$
or in one of the intervals $[x_{k},x_{k+1}]$,
and every element of $B$ lies in $[y_1,\infty)$
or in one of the intervals $[y_{k+1},x_{k}]$.
Now consider the continuous map $h$ that
has slope $\frac12$ on $(-\infty,x_1]$ and on $[y_1,\infty)$,
sends $x_k$ to $x_{k+1}$ and $y_k$ to $y_{k+1}$ for all $k$,
and is piecewise linear away from the points $x_k$, $y_k$;
it is well-defined because these intervals cover $R$,
and by looking at its behavior on each of those intervals
we can see that it has no fixed points.
On the other hand, $h$ is a contraction map 
with contraction constant $\frac12$.  Contradiction.

\bigskip

The Alternating Series Test (14) does not imply completeness:
In the field of formal Laurent series,
every series whose terms tend to zero
(whether or not they alternate in sign)
is summable,
so the Alternating Series Test holds
even though the Cut Property doesn't.

\bigskip

The Absolute Convergence Property (15) does not imply completeness:
The field of formal Laurent series has the property
that every absolutely convergent series is convergent
(and indeed the reverse is true as well!),
but it does not satisfy the Cut Property.

\bigskip

The Ratio Test (16) implies completeness by way of 
the Archimedean Property (2) and the Cut Property (3):
Note that the Ratio Test implies that
$\frac12 + \frac14 + \frac18 + \dots$ converges,
implying that $R$ is Archimedean
(the sequence of partial sums $\frac12$, $\frac34$, $\frac78$, \dots
isn't even a Cauchy sequence if there exists
an $\epsilon > 0$ that is less than $1/n$ for all $n$).
Now we make use of the important fact
(which we have avoided making use of up till now, for esthetic reasons,
but which could be used to expedite some of the preceding proofs)
that every Archimedean ordered field is isomorphic to a subfield of the reals.
(See the next paragraph for a proof.)
To show that a subfield of the reals that satisfies the Ratio Test
must contain every real number, it suffices to note
that every real number can be written as a sum
$n \pm \frac12 \pm \frac14 \pm \frac18 \pm \dots$
that satisfies the hypotheses of the Ratio Test.

\bigskip

Every Archimedean ordered field is isomorphic to a subfield of the reals:
For every $x$ in $R$, let $S_x$ be the set of elements of $\bQ$
whose counterparts in $\bQ_R$ are less than $x$,
and let $\phi(x)$ be the least upper bound of $S_x$.
The Archimedean Property can be used to show that $\phi$ is an injection,
and with some work one can verify that it is also a field homomorphism.
For more on completion of ordered fields, see \cite{refSb}.

\bigskip

The Shrinking Interval Property (17) does not imply completeness:
The field of formal Laurent series satisfies the former but not the latter.
For details, see \cite[pp.\ 212--215]{refEf}.

\bigskip

The Nested Interval Property (18) does not imply completeness:
The surreal numbers are a counterexample.
(Note however that the field of formal Laurent series
is {\it not\/} a counterexample;
although it satisfies the Shrinking Interval Property,
it does not satisfy the Nested Interval Property,
since for instance the nested closed intervals $[n,\omega/n]$
have empty intersection.
This shows that, as an ordered field property,
(18) is strictly stronger than (17).)
To verify that the surreal numbers satisfy (18),
consider a sequence of nested intervals 
$[a_1,b_1] \supseteq [a_2,b_2] \supseteq \dots$
If $a_i = b_i$ for some $i$, say $a_i = b_i = c$,
then $a_j = b_j = c$ for all $j>i$,
and $c$ lies in all the intervals.
If $a_i < b_i$ for every $i$,
then $a_i \leq a_{\max\{i,j\}} < b_{\max\{i,j\}} \leq b_j$ for all $i,j$,
so every element of $A = \{a_1, a_2, \dots\}$ 
is less than every element of $B = \{b_1, b_2, \dots\}$.
Hence there exists a surreal number
that is greater than every element of $A$
and less than every element of $B$,
and this surreal number lies in all the intervals $[a_i,b_i]$.
Thus $\bNo$ satisfies the Nested Interval Property
but being non-Archimedean does not satisfy completeness.

If one dislikes this counterexample
because the surreal numbers are a Field rather than a field,
one can instead use the field of surreal numbers
that ``are created before Day $\omega_1$'',
where $\omega_1$ is the first uncountable ordinal.
See \cite{refCo} for a discussion of the ``birthdays'' of surreal numbers.
For a self-contained explanation of a related counterexample
that predates Conway's theory of surreal numbers, see \cite{refD}.
For a counterexample arising from non-standard analysis,
see the discussion of the Cantor completeness of 
Robinson's valuation field $^\rho \mathbb R$
in \cite{refHa} and \cite{refHT}.

Because these counterexamples are abstruse,
one can find in the literature and on the web
assertions like ``The Nested Interval Property implies 
the Bolzano-Weierstrass Theorem and vice versa''.
It's easy for students to appeal to the Archimedean Property
without realizing they are doing so,
especially because concrete examples of non-Archimedean ordered fields
are unfamiliar to them.

%Mention the version of the NIP that allows arbitrary (not just countable)
%collections of intervals, and maybe the 2-Helly property too?

\bigskip

To summarize: Properties (1), (3), (4), (5), (7), (8), (9), 
(10), (12), (13), and (16) imply completeness,
while properties (2), (6), (11), (14), (15), (17), and (18) don't.
The ordered field of formal Laurent series
witnesses the fact that (11), (14), (15), and (17) don't imply completeness;
some much bigger ordered fields witness the fact 
that (6) and (18) don't imply completeness;
and every non-Archimedean ordered field witnesses the fact
that (2) doesn't imply completeness.

One of the referees asked which of the properties 
(6), (11), (14), (15), (17), and (18)
imply completeness in the presence of the Archimedean Property (2).
The answer is, All of them.
It is easy to show this in the case
of properties (11), (14), (15), (17), and (18),
using the fact that every Archimedean ordered field
is isomorphic to a subfield of the reals
(see the discussion of property (16) in section 4). 
The case of (6) is slightly more challenging.

\bigskip

Claim: Every Archimedean ordered field with the Bounded Value Property
is Dedekind complete.

Proof: Suppose not; let $R$ be a counterexample,
and let $A,B$ be a bad cut of $R$.
Let $a_n = \max A \cap 2^{-n} \bZ_R$
and $b_n = \min B \cap 2^{-n} \bZ_R$,
so that $|a_n - b_n| = 2^{-n}$.
Since $A,B$ is a bad cut,
neither of the sequences $(a_n)_{n=1}^{\infty}$ and $(b_n)_{n=1}^{\infty}$
can be eventually constant.
Let $f_n(x)$ be the continuous function that 
is 0 on $(-\infty,a_n\!-\!2^{-n}]$, 
1 on $[a_n,b_n]$, 0 on $[b_n\!+\!2^{-n},\infty)$,
and piecewise linear on $[a_n\!-\!2^{-n},a_n]$ and $[b_n,b_n\!+\!2^{-n}]$.
The interval $[a_n-2^{-n},b_n+2^{-n}]$ has length $\leq 3 \cdot 2^{-n}$,
which goes to 0 in $R$ as $n \rightarrow \infty$ since $R$ is Archimedean.
Any $c$ belonging to all the intervals $[a_n-2^{-n},b_n+2^{-n}]$
would be a cutpoint for $A,B$,
and since the cut $A,B$ has been assumed to have no cutpoint,
$\cap_{n=1}^{\infty} [a_n\!-\!2^{-n},b_n\!+\!2^{-n}]$ is empty.
It follows that for every $x$ in $R$ only finitely many 
of the intervals $[a_n\!-\!2^{-n},b_n\!+\!2^{-n}]$ contain $x$,
so $f(x) = \sum_{n=1}^{\infty} f_n(x)$ is well-defined for all $x$
(since all but finitely many of the summands vanish).
Furthermore, the function $f(x)$ is continuous,
since a finite sum of continuous functions is continuous,
and since for every $x$ we can find an $m$ and an $\epsilon>0$
such that $f_n(y) = 0$ for all $n > m$ and all $y$ with $|x-y| < \epsilon$
(so that $f$ agrees with the continuous function $\sum_{n \leq m} f_n$
on a neighborhood of $x$).
Finally note that $f$ is unbounded,
since e.g.\ $f(x) \geq n$ for all $x$ in $[a_n,b_n]$.

Alternatively, one can argue as follows:
The Archimedean Property implies countable cofinality
(specifically, $\bN_R$ is a countable unbounded set),
and an argument of Teismann \cite{refT} shows
that every ordered field with countable cofinality 
that satisfies property (6) is complete.

\bigskip

It is worth noting that for all 18 of 
the propositions listed in section 2, the answer to the question 
``Does it imply the Dedekind Completeness Property?''\ (1) 
is the same as the answer to the question 
``Does it imply the Archimedean Property?''\ (2).
A priori, one might have imagined that one or more
of properties (3) through (18)
would be strong enough to imply the Archimedean Property
yet not so strong as to be a completeness property for the reals.

\bigskip

This is not the end of the story of real analysis in reverse;
there are other theorems in analysis with which
one could play the same game.
Indeed, some readers may already be wondering
``What about the Fundamental Theorem of Calculus?''
Actually, the FTC is really two theorems, not one
(sometimes called FTC I and FTC II in textbooks).
They are not treated here
because this essay is already on the long side for a Monthly article,
and a digression into the theory of the Riemann integral
would require a whole section in itself.
Indeed, there are different ways of defining the Riemann integral
(Darboux's and Riemann's come immediately to mind),  
and while they are equivalent in the case of the real numbers,
it is possible that different definitions of the Riemann integral
that are equivalent over the reals
might turn out to be different over ordered fields in general;
thus one might obtain different varieties of FTC I and FTC II,
some of which would be completeness properties
and others of which would not.
It seemed best to leave this topic for others to explore.

An additional completeness axiom for the reals is
the ``principle of real induction'' \cite{refCl}.

\section{Some odds and ends}

\subsection{History and terminology}

It's unfortunate that the word completeness is used
to signify both Cauchy completeness and Dedekind completeness;
no doubt this doubleness of meaning has contributed to the misimpression
that the two are equivalent in the presence of the ordered field axioms.
It's therefore tempting to sidestep the ambiguity of the word ``complete''
by resurrecting Dedekind's own terminology (``{\it Stetigkeit\/}'')
and referring to the completeness property of the reals
as the {\bf continuity\/} property of the reals ---
where here we are to understand the adjective ``continuous''
not in its usual sense, as a description of a certain kind of function,
but rather as a description of a certain kind of set,
namely, the kind of set that deserves to be called a continuum.
However, it seems a bit late in the day to try to get people
to change their terminology.

It's worth pausing here to explain what Hilbert had in mind
% solicit input from a genuine historian of mathematics
when he referred to the real numbers as 
the ``complete Archimedean ordered field''.
What he meant by this
is that the real numbers can be characterized by
a property referred to earlier in this article
(after the discussion of property (16) in section 4):
every Archimedean ordered field is isomorphic
to a subfield of the real numbers.
That is, every Archimedean ordered field can be embedded
in an Archimedean ordered field that is isomorphic to the reals,
and no further extension to a larger ordered field is possible
without sacrificing the Archimedean Property.
Hilbert was saying that
the real number system is the (up to isomorphism) unique
Archimedean ordered field that is not a proper subfield of
a larger Archimedean ordered field;
vis-a-vis the ordered field axioms and the Archimedean Property,
$\bR$ is complete in the sense that nothing can be added to it.
In particular Hilbert was not asserting
any properties of $\bR$ as a metric space.

Readers interested in the original essays of Dedekind and Hilbert
may wish to read Dedekind's ``Continuity and irrational numbers'' 
and Hilbert's ``On the concept of number'', 
both of which can be found in \cite{refEw}.

\subsection{Advantages and disadvantages of the Cut Property}

The symmetry and simplicity of the Cut Property
have already been mentioned.
Another advantage is shallowness.
Although the word has a pejorative sound,
shallowness in the logical sense can be a good thing;
a proposition with too many levels of quantifiers in it
is hard for the mind to grasp.
The proposition ``$c$ is an upper bound of $S$''
(i.e., ``for all $s \in S$, $s \leq c$'')
involves a universal quantifier,
so the proposition ``$c$ is the least upper bound of $S$''
involves two levels of quantifiers,
and the proposition that for every nonempty bounded set $S$ of reals
there exists a real number $c$
such that $c$ is a least upper bound of $S$
therefore involves four levels of quantifiers.

In contrast, the assertion that $A,B$ is a cut of $R$
involves one level of quantifiers,
and the assertion that $c$ is a cutpoint of $A,B$
involves two levels of quantifiers,
so the assertion that every cut of $R$ determines a cutpoint
involves only three levels of quantifiers.

Note also that the objects 
with which the Dedekind Completeness Property is concerned ---
arbitrary nonempty bounded subsets of the reals ---
are hard to picture, 
whereas the objects with which the Cut Property is concerned ---
ways of dividing the number line into a left-set and a right-set ---
are easy to picture.  
Indeed, in the context of foundations of geometry,
it is widely acknowledged that
some version of the Cut Property is the right way to capture
what Dedekind called the continuity property of the line.

It should also be mentioned here that
the Cut Property can be viewed as
a special case of the Least Upper Bound Property,
where the set $S$ has a very special structure.
This makes the former more suitable for doing naive reverse mathematics
(since a weaker property is easier to verify)
but also makes it slightly less convenient for doing forward mathematics
(since a weaker property is harder to use).
If one starts to rewrite a real analysis textbook 
replacing every appeal to the Least Upper Bound Property
by an appeal to the Cut Property,
one quickly sees that one ends up mimicking the textbook proofs
but with extra, routinized steps
(``\dots and let $B$ be the complement of that set'')
that take up extra space on the page and add no extra insight.
So even if one wants to assign primacy to the Cut Property,
one would not want to throw away the Least Upper Bound Property;
one would introduce it as a valuable consequence of the Cut Property.

Lastly, we mention a variant of the Cut Property, 
Tarski's Axiom 3 \cite{refWc}, that drops the hypothesis that 
the union of the two sets is the whole ordered field $R$.
This stronger version of the axiom is equivalent to the one presented above.
Like the Least Upper Bound Property, Tarski's version of the Cut Property 
is superior for the purpose of constructing the theory of the reals
but less handy for the purpose of ``deconstructing'' it.

\subsection{Implications for pedagogy}

As every thoughtful teacher knows,
logical equivalence is not the same as pedagogical equivalence.
Which completeness property of the reals should we teach
to our various student audiences, assuming we teach one at all?

Here the author drops the authorial ``we''
(appropriate for statements of a mathematical and historical nature
that are, as far as the author has been able to assess, accurate)
and adopts an authorial ``I'' more appropriate to statements of opinion.

I think the reader already knows
that I am quite taken with the Cut Property as an axiom for the reals,
and will not be surprised to hear that
I would like to see more teachers of calculus,
and all teachers of real analysis,
adopt it as part of the explicit foundation of the subject.
What may come as a bigger surprise is that
I see advantages to a different completeness axiom
that has not been mentioned earlier in the article,
largely because I have not seen it stated in any textbook
(although Burns \cite{refBu} does something similar, as I describe below).

When we write $.3333\dots$, what we mean 
(or at least one thing we mean) is
``The number that lies between .3 and .4, and lies between .33 and .34,
and lies between .333 and .334, etc.''
That is, {\em a decimal expansion is
an ``address'' of a point of the number line}.
Implicit in the notation is the assumption that for every decimal expansion,
such a number {\it exists\/} and is {\it unique\/}.
These assumptions of existence and uniqueness
are part of the mathematical undermind
(the mathematical subconscious, if you prefer)
of the typical high schooler.
After all, it never occurs to a typical high school student
whether there might be more than one number 0.5,
or whether there might be no such number at all
(though balking at fractions is common
for thoughtful students at an earlier age);
% Ask Donna if this is true
so it's tempting to carry over the assumption of existence and uniqueness
when the teacher makes the transition from finite decimal to infinite decimals. 

Part of what an honors calculus teacher should do
is undermine the mathematical undermind, and convince the students
that they've been uncritically accepting
precepts that have not yet been fully justified.
The flip side of von Neumann's adage
``In mathematics you don't understand things; you just get used to them''
is that once you get used to something,
you may mistakenly come to believe you understand it!
Infinite decimals can come to seem intuitive,
on the strength of their analogy with finite decimals,
and the usefulness of infinite decimals
makes us reluctant to question the assumptions
on which they are based.
But mathematics is a liberal art,
and that means we should bring difficulties into the light
and either solve them honestly or duck them honestly.
And the way a mathematician ducks a problem honestly
is to formulate the problematic assumption
as precisely and narrowly as possible and call it an axiom.

Specifically, I would argue that one very pedagogically appropriate axiom
for the completeness of the reals
is one that our students have been implicitly relying on for years:
the Strong Nested Decimal Interval Property,
which asserts that for all infinite strings
$d_0,d_1,d_2,\dots$ of digits between 0 and 9,
there exists a unique real number in the intervals
$[.d_0, .d_0 + .1]$, $[.d_0 d_1, .d_0 d_1 + .01]$, 
$[.d_0 d_1 d_2, .d_0 d_1 d_2 + .001]$, etc.
(I call it ``Strong'' because, unlike the ordinary Nested Interval Property,
it asserts uniqueness as well as existence.)
The reader who has made it through the article thus far
should have no trouble verifying that
this is indeed a completeness property of the reals,
and one can use it to give expeditious proofs 
of some of the important theorems of the calculus,
at least in special cases.
(Example: To show that the Intermediate Value Theorem
holds for weakly increasing functions,
one can home in on the place where the function vanishes
by considering decimal approximations from both sides.)
This choice of axiom does not affect the fact
that the main theorems of the calculus have proofs
that are hard to understand for someone
who is taking calculus for the first or even the second time
and who does not have much practice in reading proofs;
indeed, I would say that the art of reading proofs
goes hand-in-hand with the art of writing them, 
and very few calculus students understand the forces at work
and the constraints that a mathematician labors under
when devising a proof.
But if we acknowledge early in the course that
the Strong Nested Decimal Interval Property
(or something like it) is an assumption that our theorems rely upon,
and stress that it cannot be proved by mere algebra,
we will be giving our students a truer picture of the subject.

Furthermore, the students will encounter infinite decimals 
near the end of the two-semester course when infinite series are considered;
now an expression like $.3333\dots$
means $3 \times 10^{-1} + 3 \times 10^{-2} + 3 \times 10^{-3}
+ 3 \times 10^{-4} + \dots$.
(Burns \cite{refBu} adopts as his completeness axioms
the Archimedean Property plus the assertion that every decimal converges.)
The double meaning of infinite decimals hides a nontrivial theorem:
Every infinite decimal, construed as an infinite series,
converges to a limit, specifically,
the unique number that lies in all the associated nested decimal intervals.
We do not need to prove this assertion
to give our students the knowledge that this assertion has nontrivial content;
we can lead them to see that the calculus gives them, for the first time,
an honest way of seeing why $.9999\dots$ is the same number as $1.0000\dots$
(a fact that they may have learned to parrot
but probably don't feel entirely comfortable with).

Students should also be led to see that
the question ``But how do we know that the square root of two really exists,
if we can't write down all its digits or give a pattern for them?''\ is
a fairly intelligent question.
In what sense do we know that such a number exists?
We can construct the square root of two as the length of the diagonal 
of a square of side-length one,
but that trick won't work if we change the question to
``How do we know that the {\it cube\/} root of two really exists?''

To those who would be inclined to show students
a construction of the real numbers
(via Dedekind cuts and Cauchy sequences),
I would argue that a student's first exposure to rigorous calculus
should focus on other things. 
It takes a good deal of mathematical sophistication
to even appreciate why someone would want
to prove that the theory of the real numbers is consistent,
and even more sophistication to appreciate
why we can do so by making a ``model'' of the theory.
Most students enter our classrooms
with two workable models of the real numbers,
one geometrical (the number line)
and one algebraic (the set of infinite decimals).
Instead of giving them a third picture of the reals,
it seems better to clarify the pictures that they already have,
and to assert the link between them.

In fact, I think that for pedagogical purposes,
it's best to present both the Cut Property
and the Strong Nested Decimal Interval Property,
reflecting the two main ways students think about real numbers.
And it's also a good idea to mention
that despite their very different appearances,
the two axioms are deducible from one another,
even though neither is derivable
from the principles of high school mathematics.
This will give the students a foretaste of a refreshing phenomenon 
that they will encounter over and over 
if they continue their mathematical education:
two mathematical journeys that take off in quite different directions
can unexpectedly lead to the same place.

\paragraph{Acknowledgments.}  
Thanks to Matt Baker, Mark Bennet, Robin Chapman, Pete Clark, Ricky Demer, 
Ali Enayat, James Hall, Lionel Levine, George Lowther, David Speyer, and 
other contributors to {\tt MathOverflow} (\url{http://www.mathoverflow.net})
for helpful comments; thanks to Wilfried Sieg for his historical insights;
thanks to the referees for their numerous suggestions of ways to make this 
article better; and special thanks to John Conway for helpful conversations. 
Thanks also to my honors calculus students, 2006--2012, whose diligence and 
intellectual curiosity led me to become interested in the foundations of the 
subject.

\raggedright

\noindent\textbf{James G. Propp} 
received his A.B.\ from Harvard College in 1982,
where he took both honors freshman calculus and real analysis
with Andrew Gleason.
He received his Ph.D.\ from the University of California--Berkeley in 1987.  
His research is mostly in the areas of combinatorics and probability,
and he currently teaches at the University of Massachusetts--Lowell.

\noindent\textit{Department of Mathematics,
University of Massachusetts--Lowell, Lowell, MA 01854\\
{\rm Website:} {\tt http://jamespropp.org}}
% I don't want to publish my email address where robots can find it; 
% actual human beings who want it can go to my web-page and easily 
% find out how to contact me by email.


\begin{thebibliography}{25}

\bibitem{refAS} C.\ Apelian, S.\ Surace, \emph{Real and Complex Analysis}.
Chapman \& Hall / CRC Pure and Applied Mathematics, 2009.

\bibitem{refBr} A.\ Browder, \emph{Mathematical Analysis},
Undergraduate Texts in Mathematics.  Springer-Verlag, New York, 1996.

\bibitem{refBu} R.\ P.\ Burns, \emph{Numbers and Functions: Steps to Analysis}.
Cambridge University Press, 1992 and 2000.

\bibitem{refCl} P.\ L.\ Clark, The instructor's guide to real induction, 2012,
available at \url{http://arxiv.org/abs/1208.0973}.

\bibitem{refCo} J.\ H.\ Conway, \emph{On Numbers and Games}. 
Academic Press, London, New York, San Francisco, 1976.

\bibitem{refD} R.\ O.\ Davies, solution to advanced problem 5112,
\emph{Amer. Math. Monthly} \textbf{72} (1965) 85--87.

\bibitem{refEf} N.\ V.\ Efimov, \emph{Higher Geometry}.  Mir, Moscow, 1980.

\bibitem{refEn} A.\ Enayat,
\url{http://mathoverflow.net/questions/71432/ordered-fields-with-the-bounded-value-property}.

\bibitem{refEw} W.\ Ewald,
\emph{From Kant to Hilbert: A Source Book in the Foundations of Mathematics},
vol.\ II.  Oxford University Press, 1996.

\bibitem{refG} M.\ J.\ Greenberg,
Old and new results in the foundations of elementary 
plane Euclidean and non-Euclidean geometry, 
\emph{Amer. Math. Monthly} \textbf{117} (2010) 198--219.

\bibitem{refHa} J.\ F.\ Hall, Completeness of ordered fields, 2010,
available at \url{http://arxiv.org/abs/1101.5652}.

\bibitem{refHT} J.\ F.\ Hall, T.\ D.\ Todorov,
Completeness of the Leibniz field and rigorousness of infinitesimal calculus, 
2011, available at \url{http://arxiv.org/abs/1109.2098}.

\bibitem{refHb} P.\ R.\ Halmos, \emph{Naive Set Theory}.
Undergraduate Texts in Mathematics, Springer-Verlag, 1960.

\bibitem{refJS}
H.\ J.\ Keisler, J.\ H.\ Schmerl,
Making the hyperreal line both saturated and complete. 
\emph{J. Symbolic Logic} \textbf{56} (1991) 1016--1025. 

\bibitem{refKa} T.\ W.\ Korner, \emph{A Companion to Analysis: 
A Second First and First Second Course in Analysis}.
Graduate Studies in Mathematics, American Mathematical Society, 2004.

\bibitem{refKb} S.\ G.\ Krantz, \emph{An Episodic History of Mathematics:
Mathematical Culture through Problem Solving}, 2006, available at
\url{http://www.math.wustl.edu/~sk/books/newhist.pdf}.

\bibitem{refL} G.\ Lowther,
\url{http://mathoverflow.net/questions/65874/converse-to-banachs-fixed-point-theorem-for-ordered-fields}.

\bibitem{refDM} D.\ Marker, \emph{Model Theory: An Introduction},
Graduate Texts in Mathematics, vol.\ 217, Springer-Verlag, 2002.

\bibitem{refGM} G.\ E.\ Martin, \emph{The Foundations of Geometry
and the Non-Euclidean Plane},
Undergraduate Texts in Mathematics, Springer-Verlag, 1975.

\bibitem{refMe} E.\ Mendelson,
\emph{Number Systems and the Foundations of Analysis}.  Dover, 1973.

% \bibitem{refMR} R.\ M.\ F.\ Moss and G.\ T.\ Roberts, A creeping lemma,
% \emph{Amer.\ Math.\ Monthly} {\bf 75} (1968), 649--652.

\bibitem{refP} J.\ Propp,
A frosh-friendly completeness axiom for the reals, preprint (2009),
available at \url{http://jamespropp.org/cut.pdf}.

\bibitem{refSa} J.\ H.\ Schmerl,
Models of Peano arithmetic and a question of Sikorski on ordered fields.
\emph{Israel J. Math.} \textbf{50} (1985) 145--159.

\bibitem{refSb} D.\ Scott, On completing ordered fields.
In: \emph{Applications of Model Theory to Algebra, Analysis, and Probability},
Holt, Rinehart, and Winston, New York, 1969.  274--278.

\bibitem{refSc} R.\ Sikorski, On an ordered algebraic field. 
\emph{Towarzytwo Nankowe Warzawskie} \textbf{41} (1948) 69--96.

\bibitem{refSd} S.\ G.\ Simpson, \emph{Subsystems of second order arithmetic}, 
Second edition. Perspectives in Logic, Cambridge University Press, Cambridge; 
Association for Symbolic Logic, Poughkeepsie, NY, 2009.

\bibitem{refT} H.\ Teismann,
Towards a (more) complete list of completeness axioms,
\emph{Amer. Math. Monthly} \textbf{120} (2013) 99--114.

\bibitem{refWa} 
Wikipedia contributors, 
Construction of the real numbers,
\emph{Wikipedia, The Free Encyclopedia}, available at
\url{http://en.wikipedia.org/wiki/Construction_of_the_real_numbers};
retrieved March 30, 2012.

\bibitem{refWaa}
Wikipedia contributors,
Dedekind-MacNeille completion,
\emph{Wikipedia, The Free Encyclopedia}, available at
\url{http://en.wikipedia.org/wiki/Dedekind-MacNeille_completion};
retrieved July 20, 2012.

\bibitem{refWb} 
Wikipedia contributors, 
Infinitesimal,
\emph{Wikipedia, The Free Encyclopedia}, available at
\url{http://en.wikipedia.org/wiki/Infinitesimal};
retrieved March 30, 2012.

\bibitem{refWc}
Wikipedia contributors, 
Tarski's axiomatization of the reals,
\emph{Wikipedia, The Free Encyclopedia}, available at
\url{http://en.wikipedia.org/wiki/Tarski's_axiomatization_of_the_reals};
retrieved March 30, 2012.

\end{thebibliography}
\end{document}